\documentclass[12pt]{article}
\usepackage{amsmath}
\usepackage{amsthm}
\usepackage{amsfonts}
\usepackage{amssymb}
\usepackage{esint}
\usepackage{eucal}
\usepackage{mathrsfs}
\usepackage{graphicx,graphics}
\numberwithin{equation}{section}
\usepackage{graphicx,graphics}
\usepackage{pdfsync}
\usepackage{color}
\def \dis {\displaystyle}
\def \tex {\textstyle}

\def \confai {-\kern -.5em\rightharpoonup}
\def \cqfd {\hfill$\Box$}

\def \al {\alpha}
\def \be {\beta}
\def \ga {\gamma}

\def \De {\Delta}
\def \ep {\varepsilon}

\def \Om {\Omega}
\def \la {\lambda}

\def \ph {\varphi}

\def \si {\sigma}

\def \ZZ {\mathbb Z}
\def \QQ {\mathbb Q}
\def \RR {\mathbb R}

\def \beq {\begin{equation}}
\def \eeq {\end{equation}}
\def \ba {\begin{array}}
\def \ea {\end{array}}
\def \bs {\bigskip}

\def \ecart {\noalign{\medskip}}
\newtheorem{Thm}{Theorem}[section]
\newtheorem{Cor}[Thm]{Corollary}
\newtheorem{Pro}[Thm]{Proposition}

\newtheorem{Adef}[Thm]{Definition}

\newtheorem{Arem}[Thm]{Remark}
\newenvironment{Rem}{\begin{Arem}\rm}{\end{Arem}}
\newtheorem{Aexa}[Thm]{Example}
\newenvironment{Exa}{\begin{Aexa}\rm}{\end{Aexa}}
\newtheorem{Anot}[Thm]{Notation}

\def \refe #1.{\eqref{#1})}
\def \reff #1.{figure~\ref{#1}}
\def \refs #1.{Section~\ref{#1}}
\def \refss #1.{Subsection~\ref{#1}}
\def \refD #1.{Definition~\ref{#1}}
\def \refT #1.{Theorem~\ref{#1}}
\def \refL #1.{Lemma~\ref{#1}}
\def \refC #1.{Corollary~\ref{#1}}
\def \refP #1.{Proposition~\ref{#1}}
\def \refPt #1.{Properties~\ref{#1}}
\def \refR #1.{Remark~\ref{#1}}
\def \refE #1.{Example~\ref{#1}}
\def \refN #1.{Notation~\ref{#1}}
\newcounter{marnote}

\usepackage[colorlinks=true, linkcolor=blue]{hyperref} 

\title{Isotropic realizability of fields and reconstruction of invariant measures under positivity properties.
Asymptotics of the flow by a non-ergodic approach}
\author{Marc Briane
\\
\normalsize Univ Rennes, INSA Rennes,  CNRS, IRMAR - UMR 6625, F-35000 Rennes, France
\\
\normalsize mbriane@insa-rennes.fr
}
\begin{document}
\maketitle
\begin{abstract}
The paper is devoted to the isotropic realizability of a regular gradient field $\nabla u$ or a more general vector field $b$, namely the existence of a continuous positive function $\sigma$ such that $\sigma b$ is divergence free in $\RR^d$ or in an open set of $\RR^d$. First, we prove that under some suitable positivity condition satisfied by $\nabla u$, the isotropic realizability of $\nabla u$ holds either in~$\RR^d$ if $\nabla u$ does not vanish, or in the open sets $\{c_j\!<\!u\!<\!c_{j+1}\}$ if the $c_j$ are the critical values of~$u$ (including $\inf_{\RR^d}u$ and $\sup_{\RR^d}u$) which are assumed to be in finite number. It turns out that this positivity condition is not sufficient to ensure the existence of a continuous positive invariant measure $\sigma$ on the torus when $\nabla u$ is periodic. Then, we establish a new criterium of the existence of an invariant measure for the flow associated with a regular periodic vector field $b$, which is based on the equality $b\cdot\nabla v=1$ in $\RR^d$. We show that this gradient invertibility is not related to the classical ergodic assumption, but it actually appears as an alternative to get the asymptotics of the flow.
\end{abstract}
\par\bs\noindent
\section{Introduction}
In this paper we study the problem of the isotropic realizability of a vector field $b\in C^1(\RR^d)^d$, namely the existence of a  positive function $\sigma\in C^0(\RR^d)$ solution to the equation
\beq\label{bsi}
{\rm div}\,(\sigma b)=0\;\;\mbox{in }\RR^d,
\eeq
or in an open subset of $\RR^d$. When the vector field $b$ is periodic with respect to $Y_d:=[0,1)^d$, {\em i.e.}
\[
\forall\,\kappa\in\ZZ^d,\ \forall\,x\in\RR^d,\quad b(x+\kappa)=b(x),
\]
the problem of the isotropic realizability in the torus $\RR^d/\ZZ^d$, namely the existence of a positive $Y_d$-periodic function $\sigma\in C^0(\RR^d)$ solution to \eqref{bsi}, is also addressed.
\par
In the case where $b=\nabla u$ is a gradient field, the reconstruction of a positive $\sigma$ has been first done in \cite{BoVa} and a rigorous way in \cite{Ric} assuming that $\nabla u$ never vanishes and using the method of characteristics. Alternatively, when the potential $u$ satisfies a prescribed boundary condition on a bounded smooth domain of $\RR^2$ with a finite number of critical points, a conductivity $\sigma$ has been derived in \cite{Ale} thanks to an approximation procedure adding a vanishing viscosity term. More recently, the isotropic realizability of a non-vanishing gradient, {\em i.e.}
\beq\label{Du>0}
\inf_{\RR^d}|\nabla u|>0,
\eeq
has been revisited in \cite{BMT} both in the space $\RR^d$ and in the torus $\RR^d/\ZZ^d$ using specifically the flow
\beq\label{bflow}
\left\{\ba{ll}
X'(t,x)=b\big(X(t,x)\big), & t\in\RR
\\ \ecart
X(0,x)=x\in\RR^d,
\ea\right.
\eeq
with $b=\nabla u$.
In particular, it was proved that the isotropic realizability in the torus is actually stronger that the realizability in the space.
Furthermore, again using the flow \eqref{bflow} we showed in~\cite[Theorem~4.1]{Bri1} that in any dimension the presence of critical points for the potential $u$ may be an obstacle to the (even local) existence of a conductivity $\sigma$ solution to \eqref{bsi} with $b=\nabla u$. The case of non-regular gradients has been also investigated in \cite{Bri2}.
\par\bs
Beyond the negative results of \cite{Bri1} when the non-vanishing condition \eqref{Du>0} does not hold, we thus need extra conditions on the gradient $\nabla u$ to ensure its isotropic realizability in the whole space $\RR^d$ or at least in a subset of $\RR^d$. First, we prove (see Theorem~\ref{thm.irg}) that if $u\in C^2(\RR^d)$ has either all its non-negative partial derivatives or all its non-positive partial derivatives the ratios of which are controlled from above and below (see more precisely condition~\eqref{disdku} below), and if $u$ has exactly $n$ critical values:
\[
c_0:=\inf_{\RR^d}u<c_1=u(\xi^1)<\cdots<c_n=u(\xi^n)<c_{n+1}:=\sup_{\RR^d}u,\quad\mbox{with }\nabla u(\xi^j)=0,
\]
possibly with $\inf_{\RR^d}|\nabla u|=0$ (so that $\nabla u$ may vanish at infinity), then $\nabla u$ is isotropically realizable:
\begin{itemize}
\item either in $\RR^d$ when $u$ has no critical point,
\item or in the $(n\!+\!1)$ open sets $\{c_j<u<c_{j+1}\}$ for $j=0,\dots,n$.
\end{itemize}
Then, we extend this result to the isotropic realizability of a vector field $b\in C^1(\RR^d)^d$. Assuming the existence of an open interval $I\subset\RR$ and a function $u\in C^1(\RR^d)$ such that for any $x$ in the inverse image $\{u\in I\}$, the function $u(X(\cdot,x))$
is increasing and its range contains $I$, we show (see Theorem~\ref{thm.irb}) the isotropic realizability of $b$ in the open set $\{u\in I\}$.
\par\bs
The isotropic realizability in the torus $\RR^d/\ZZ^d$ of a $Y_d$-periodic vector field $b\in C^1(\RR^d)^d$ is more intricate.
In this case by the uniqueness of the Cauchy-Lipschitz theorem the flow $X$ solution to \eqref{bflow} satisfies
\[
\forall\,\kappa\in\ZZ^d,\ \forall\,(t,x)\in\RR\times\RR^d,\quad X(t,x+\kappa)=X(t,x)+\kappa,
\]
so that the image of $X(t,x)$ by the canonical surjection $\Pi:\RR^d\to\RR^d/\ZZ^d$ is independent of any representative $x$ in the class $\Pi(x)$.
Hence, equation \eqref{bflow} is well posed in the torus. Then, by virtue of Liouville theorem (see, {\em e.g.}, \cite[Chap. 2, Theorem 1]{CFS}) the isotropic realizability \eqref{bsi} in the torus is equivalent to the existence of a positive $Y_d$-periodic function $\sigma\in C^0(\RR^d)$ which is called an {\em invariant measure} for the flow $X$, such that for any $Y_d$-periodic function $\ph\in C^1(\RR^d)$,
\beq\label{imsi}
\forall\,t\in\RR,\quad\int_{Y_d}\ph(X(t,x))\,\sigma(x)\,dx=\int_{Y_d}\ph(x)\,\sigma(x)\,dx.
\eeq
Furthermore, the natural extension to any vector field $b$ of the condition \eqref{disdku} relating to a gradient field, is the boundedness from below by a positive constant of the coefficients $b_k$ for $k=1,\dots, d$, either the coefficients $-\,b_k$. However, it turns out that this boundedness condition is not sufficient to get the existence of an invariant measure as show Proposition~\ref{pro.neim} and Example~\ref{exa.noim}.
\par
Actually, assuming that the $Y_d$-periodic vector field $b\in C^1(\RR^d)^d$ satisfies the {\em gradient invertibility}
\beq\label{ibDv1}
b\cdot\nabla v_1=1\;\;\mbox{in }\RR^d,
\eeq
for some $Y_d$-periodic gradient $\nabla v_1\in C^0(\RR^d)^d$, we prove (see Theorem~\ref{thm.eim}) that the existence of an invariant measure \eqref{imsi} for the flow $X$ is equivalent to the existence of a vector $\xi\in\RR^d$ and $d$ linearly independent $Y_d$-periodic gradients $\nabla w_k\in C^0(\RR^d)^d$ such that
\beq\label{bDwk}
b\cdot\nabla w_k=\xi_k\;\;\mbox{in }\RR^d,\quad\mbox{for }k=1,\dots,d.
\eeq
At this point, we need to replace in dimension $d\geq 3$ the equation \eqref{bsi} by the more restrictive condition that $b$ is proportional to a cross product of $(d\!-\!1)$ gradients. As a by-product, under condition \eqref{ibDv1} the former equivalence shows (see Corollary~\ref{cor.asyX}) that the existence of an invariant measure for $X$ implies the asymptotics
\[
\lim_{|t|\to\infty}{X(t,x)\over t}=\xi\quad\mbox{for any }x\in\RR^d,
\]
and not only almost everywhere in $\RR^d$ as obtained by the Birkhoff ergodic theorem. Surprisingly, although the limit $\xi$ is constant, it appears (see Example~\ref{exa.nonergo}) that the flow $X$ is not in general ergodic in dimension $d\geq 2$. Indeed, we may construct a non-constant $Y_d$-periodic function which is invariant by the flow $X$. Therefore, it seems that the gradient invertibility \eqref{ibDv1} can be regarded as a substitute for the classical ergodic assumption (see Remark~\ref{rem.asyX}). This allows us to recover some of the two-dimensional ergodicity results of \cite{Tas,Pei} by a new and non-ergodic approach, and to extend partially them to higher dimension.
As a natural extension of Corollary~\ref{cor.asyX} the homogenization of a linear transport equation with oscillating coefficients (see Corollary~\ref{cor.htranspeq}) is derived by the non-ergodic approach. Condition~\eqref{ibDv1} in any dimension still plays the same role as the irrationality of the so-called {\em rotation number} (see Remark~\ref{rem.asyX}~3.) in the two-dimensional homogenization results of \cite{Bre,HoXi,Tas} which are based on the ergodicity of the flow.
\subsubsection*{Notations}
\begin{itemize}
\item $\left(e_1,\dots,e_d\right)$ denotes the canonical basis of $\RR^d$.
\item $\cdot$ denotes the scalar product in $\RR^d$.
\item $I_d$ denotes the unit matrix of $\RR^{d\times d}$, and $R_\perp$ denotes the clockwise $90^\circ$ rotation matrix in $\RR^{2\times 2}$.
\item For $M\in\RR^{d\times d}$, $M^T$ denotes the transpose of $M$.
\item $Y_d:=\left[0,1\right)^d$, and $\langle f\rangle$ denotes the average-value of a function $f\in L^1(Y_d)$.
\item $|A|$ denotes the Lebesgue measure of a measurable subset $A$ of $\RR^d$.
\item $C^k_\sharp(Y_d)$ denotes the space of the $Y_d$-periodic functions of class $C^k$ in $\RR^d$.
\item $L^p_\sharp(Y_d)$, $p\geq 1$, denotes the space of the $Y_d$-periodic functions in $L^p_{\rm loc}(\RR^d)$, and $H^1_\sharp(Y_d)$ denotes the space of the functions $\ph\in L^2_\sharp(Y_d)$ such that $\nabla\ph\in L^2_\sharp(Y_d)^d$.
\item For any open set $\Om$ of $\RR^d$, $C^\infty_c(\Om)$ denotes the space of the smooth functions with compact support in $\Om$.
\item For $u\in L^1_{\rm loc}(\RR^d)$ and $U=(U_j)_{1\leq j\leq d}\in L^1_{\rm loc}(\RR^d)^d$,
\beq
\nabla u:=\left(\partial_{x_1},\dots,\partial_{x_d}\right)\quad\mbox{and}\quad DU:=\big[\partial_{x_i}U_j\big]_{1\leq i,j\leq d}.
\eeq
\item For $\xi^1_1,\dots,\xi^d$ in $\RR^d$, the cross product $\xi^2\times\cdots\times \xi^d$ is defined by
\beq\label{cropro}
\xi^1\cdot\left(\xi^2\times\dots\times \xi^d\right)=\det\left(\xi^1,\xi^2,\dots,\xi^d\right)\quad\mbox{for }\xi^1\in\RR^d,
\eeq
where $\det$ is the determinant with respect to the canonical basis $(e_1,\dots,e_d)$, or equivalently, the $k^{\rm th}$ coordinate of the cross product is given by the $(d-1)\times(d-1)$ determinant
\beq\label{crossprod}
\left(\xi^2\times\cdots\times \xi^d\right)\cdot e_k=(-1)^{k+1}\left|\,\begin{smallmatrix}
\xi^2_1 & \cdots & \xi^d_1
\\
\vdots & \qquad\ddots\qquad & \vdots
\\
\\
\xi^2_{k-1} & \cdots & \xi^d_{k-1}
\\
\xi^2_{k+1} & \cdots & \xi^d_{k+1}
\\
\vdots & \qquad\ddots\qquad & \vdots
\\
\\
\xi^2_d & \cdots & \xi^d_d
\end{smallmatrix}\,\right|.
\eeq
\end{itemize}
\section{Isotropic realizability of a vector field in $\RR^d$ under positivity properties}
\subsection{Isotropic realizability of a gradient in $\RR^d$}
Let $u\in C^1(\RR^d)$. In this section we assume that the gradient field $b=\nabla u$ has the following positivity properties:
\beq\label{dku>0}
\forall\,k\in\{1,\dots,d\},\;\;\partial_{x_k}u\geq 0\mbox{ in }\RR^d\quad\mbox{or}\quad \forall\,k\in\{1,\dots,d\},\;\;\partial_{x_k}u\leq 0\mbox{ in }\RR^d,
\eeq
and there exist positive fonctions $\al_k,\beta_k\in C^0(\RR)$ with
\beq\label{alkbek}
\int_{0}^{\pm\infty} \al_k(t)\,dt=\int_{0}^{\pm\infty} \be_k(t)\,dt=\pm\infty,
\eeq
such that for any $x\in\RR^d$, up to renumber the coordinates $x_k$,
\beq\label{disdku}
\forall\,k\in\{1,\dots,d-1\},\quad {\al_k(x_k)\over\al_{k+1}(x_{k+1})}\,\big|\partial_{x_k}u(x)\big|\leq\big|\partial_{x_{k+1}}u(x)\big|\leq
\big|\partial_{x_k}u(x)\big|\,{\be_k(x_k)\over\be_{k+1}(x_{k+1})}.
\eeq
Note that in \eqref{disdku} the partial derivatives of $u$ may vanish but the ratios between two consecutive partial derivatives are controlled.
\par
We have the following result.
\begin{Thm}\label{thm.irg}
Let $u\in C^2(\RR^d)$ be a function satisfying conditions \eqref{dku>0} and \eqref{disdku}.
\par\smallskip\noindent
$i)$ Assume that $u$ has no critical point in $\RR^d$, {\em i.e.} $\nabla u$ does not vanish in $\RR^d$. Then, $\nabla u$ is isotropically realizable in $\RR^d$ with a positive function $\sigma\in C^1(\RR^d)$.
\par\medskip\noindent
$ii)$ Assume that $u$ has a unique critical point $x^0$, {\em i.e.} $\nabla u(x^0)=0$ and $\nabla u$ does not vanish in $\RR^d\setminus\{x^0\}$.
Then, $\nabla u$ is isotropically realizable with a positive $C^1$-function $\sigma$ in the open sets $\{u>u(x^0)\}$ and $\{u<u(x^0)\}$.
\par\medskip\noindent
$iii)$ More generally, assume that there exists a positive integer $n$ such that
\beq\label{cj}
\ba{lc}
& u\left(\big\{x\in\RR^d:\nabla u(x)=0\big\}\right)=\big\{c_1,\dots,c_n\big\}
\\ \ecart
\mbox{with} & \dis \inf_{\RR^d}u=:c_0<c_1<\cdots <c_n<c_{n+1}:=\sup_{\RR^d}u.
\ea
\eeq
Then, $\nabla u$ is isotropically realizable with a positive $C^1$-function $\sigma$ in the sets $\{c_j<u<c_{j+1}\}$ for $j=0,\dots,n$.
\end{Thm}
\begin{Exa}\label{exa.thmirg}
\hfill\par\smallskip\noindent
1. Let $u:\RR^2\to\RR$ be the function defined by
\[
u(x):=\arctan(x_1)+\arctan(x_2)\quad\mbox{for }x=(x_1,x_2)\in\RR^2.
\]
We have
\[
\forall\,x\in\RR^2,\quad{\partial_{x_2}u(x)\over\partial_{x_1}u(x)}={x_1^2+1\over x_2^2+1}
\]
such that condition \eqref{disdku} holds true with $\al_1(t)=\al_2(t)=\beta_1(t)=\beta_2(t)=t^2+1$.
\\
Therefore, $\nabla u$ is isotropically realizable in $\RR^2$ while $\inf_{\RR^2}|\nabla u|=0$.
\par\medskip\noindent
2. Let $u:\RR^3\to\RR$ be the function defined by
\[
u(x):=x_1^3+x_2^3+x_3^3+x_1^2x_2+x_1x_2^2+x_1^2x_3+x_1x_3^2+x_2^2x_3+x_2x_3^2\quad\mbox{for }x=(x_1,x_2,x_3)\in\RR^3.
\]
We have
\[
\forall\,x\in\RR^3,\quad
\left\{\ba{l}
\partial_{x_1}u(x)=3x_1^2+x_2^2+x_3^2+2x_1x_2+2x_1x_3
\\ \ecart
\partial_{x_2}u(x)=x_1^2+3x_2^2+x_3^2+2x_1x_2+2x_2x_3
\\ \ecart
\partial_{x_3}u(x)=x_1^2+x_2^2+3x_3^2+2x_1x_3+2x_2x_3.
\ea\right.
\]
The partial derivatives of $u$ thus turn to be $3$ quadratic forms on $\RR^3$ associated with $3$ symmetric matrices of $\RR^{3\times 3}$ the eigenvalues of which are
$0<2-\sqrt{3}<1<2+\sqrt{3}$.
Hence, the function $u$ has $(0,0,0)$ as unique critical point.
Moreover, we deduce that for any $x\in\RR^3\setminus\{(0,0,0)\}$,
\[
{2-\sqrt{3}\over 2+\sqrt{3}}\leq
\left\{\ba{l}
\dis {\partial_{x_2}u(x)\over\partial_{x_1}u(x)}={x_1^2+3x_2^2+x_3^2+2x_1x_2+2x_2x_3\over 3x_1^2+x_2^2+x_3^2+2x_1x_2+2x_1x_3}
\\ \ecart
\dis {\partial_{x_3}u(x)\over\partial_{x_2}u(x)}={x_1^2+3x_2^2+3x_3^2+2x_1x_3+2x_2x_3\over x_1^2+3x_2^2+x_3^2+2x_1x_2+2x_2x_3}
\ea\right\}
\leq {2+\sqrt{3}\over 2-\sqrt{3}},
\]
such that condition \eqref{disdku} holds true with constant functions $\al_1,\al_2,\al_3,\beta_1,\beta_2,\be_3$.
\par
Therefore, $\nabla u$ is isotropically realizable by a positive continuous function in the open sets $\{u>0\}$ and $\{u<0\}$.
\par\medskip\noindent
3. Let $f\in C^3(\RR)$ be an increasing function such that
\[
\{x\in\RR:f'(x)=0\}=\{0,1\}.
\]
Define the function $u\in C^3(\RR^3)$ by
\[
u(x):={f(x_1\!+\!x_2\!+\!x_3)\!+\!f(x_1\!+\!2x_2\!+\!x_3)\!+\!f(x_1\!+\!x_2\!+\!3x_3)\!+\!f(x_1\!+\!4x_2\!+\!5x_3)\over 4}\,\quad\mbox{for }x\in\RR^3.
\]
Due to the non-negativity of $f'$ it is easy to check that $\nabla u(x)=0$
\[
\ba{l}
\Leftrightarrow f'(x_1\!+\!x_2\!+\!x_3)=f'(x_1\!+\!2x_2\!+\!x_3)=f'(x_1\!+\!x_2\!+\!3x_3)=f'(x_1\!+\!4x_2\!+\!5x_3)=0
\\ \ecart
\Leftrightarrow x_1\!+\!x_2\!+\!x_3,\ x_1\!+\!2x_2\!+\!x_3,\ x_1\!+\!x_2\!+\!3x_3,\ x_1\!+\!4x_2\!+\!5x_3\in\{0,1\}
\\ \ecart 
\Leftrightarrow x_1\in\{0,1\},\ x_2=x_3=0.
\ea
\]
It follows that $(0,0,0)$ and $(1,0,0)$ are the only critical points of $u$ with $u(0,0,0)=f(0)$ and $u(1,0,0)=f(1)>f(0)$.
Moreover, we have for any $x\in\RR^3\setminus\{(0,0,0),(1,0,0)\}$,
\[
\ba{c}
\dis 1\leq {\partial_{x_2}u(x)\over\partial_{x_1}u(x)}={f'(x_1\!+\!x_2\!+\!x_3)\!+\!2f'(x_1\!+\!2x_2\!+\!x_3)\!+\!f'(x_1\!+\!x_2\!+\!3x_3)+4f'(x_1\!+\!4x_2\!+\!5x_3)\over f'(x_1\!+\!x_2\!+\!x_3)\!+\!f'(x_1\!+\!2x_2\!+\!x_3)\!+\!f'(x_1\!+\!x_2\!+\!3x_3)\!+\!f'(x_1\!+\!4x_2\!+\!5x_3)}\leq 4
\\ \ecart
\dis {1\over 4}\leq {\partial_{x_3}u(x)\over\partial_{x_2}u(x)}={f'(x_1\!+\!x_2\!+\!x_3)\!+\!f'(x_1\!+\!2x_2\!+\!x_3)\!+\!3f'(x_1\!+\!x_2\!+\!3x_3)\!+\!5f'(x_1\!+\!4x_2\!+\!5x_3)\over f'(x_1\!+\!x_2\!+\!x_3)\!+\!2f'(x_1\!+\!2x_2\!+\!x_3)\!+\!f'(x_1\!+\!x_2\!+\!3x_3)\!+\!4f'(x_1\!+\!4x_2\!+\!5x_3)}\leq 5,
\ea
\]
such that condition \eqref{disdku} holds true with constant functions $\al_1,\al_2,\al_3,\beta_1,\beta_2,\be_3$.
\par
Therefore, $\nabla u$ is isotropically realizable by a positive continuous function in the open sets $\{u<f(0)\}$, $\{f(0)<u<f(1)\}$, $\{u>f(1)\}$.
\end{Exa}
\par\noindent
{\bf Proof of Theorem~\ref{thm.irg}.} Let $u\in C^2(\RR^d)$.
\par\smallskip\noindent
{\it Proof of $i)$.} Fix $x\in\RR^d$. Let $0\in (\tau_-,\tau_+)$ be the maximal interval on which the gradient flow $X(\cdot,x)$ is solution to equation \eqref{bflow} with $b=\nabla u$.
The times $\tau_-$ and $\tau_+$ do depend on $x$, but their dependence is omitted for the sake of simplicity.
Define the function $f$
\beq\label{f}
f(t):=u(X(t,x))\quad\mbox{for }t\in (\tau_-,\tau_+).
\eeq
First, let us prove that the range of $f$ agrees with the interval $(\inf_{\RR^d}u,\sup_{\RR^d}u)$, {\em i.e.}
\beq\label{rfu}
\big\{f(t):t\in (\tau_-,\tau_+)\big\}=\tex{(\inf_{\RR^d}u,\sup_{\RR^d}u)}.
\eeq
In \cite{BMT} it is immediate that the range of $f$ is $\RR$, since the derivative $f'=|\nabla u(X(\cdot,x))|^2$ is defined over the whole interval $\RR$ and is bounded from below by a positive constant. Here, the flow $X(\cdot,x)$ is only defined on the interval $(\tau_-,\tau_+)$, and we may have $\inf_{\RR^d}|\nabla u|=0$. For the sake of simplicity we write $X(t)$ in place of $X(t,x)$ in the sequel.
\par
Assume by contradiction that the flow $X(t)$ is bounded in the neighborhood of $\tau_+$. Then, $\tau_+=\infty$, otherwise $X'(t)$ is bounded in the neighborhood of $\tau_+$ and the flow $X(t)$ could be extended beyond $\tau_+$ (see, {\em e.g.}, \cite[Section~17.4]{HSD}). Then, the derivative $f'(t)=|\nabla u(X(t))|^2$ is bounded from below by a positive constant in the neighborhood of $\infty$, which implies that $f(t)=u(X(t))$ tends to $\infty$ as $t\to\infty$, a contradiction.
Therefore, there exists an increasing sequence $t_n\geq 0$ which tends to $\tau_+$ such that $|X(t_n)|$ tends to~$\infty$ as $n\to\infty$.
\par
From now on, we assume that all the partial derivatives of $u$ are non-negative. The non-positivity case of condition~\eqref{dku>0} is quite similar.
Denote by $A_k$ (respectively $B_k$) a primitive of the function $\al_k$ (respectively~$\beta_k$) in condition \eqref{disdku}. We have for any $k\in\{1,\dots,d-1\}$,
\beq\label{FkGk}
\left\{\ba{ll}
A_k\big(X_k(t_n)\big)-A_{k+1}\big(X_{k+1}(t_n)\big) & \kern -.25cm \leq A_k(x_k)-A_{k+1}\big(x_{k+1})
\\ \ecart
B_k\big(X_k(t_n)\big)-B_{k+1}\big(X_{k+1}(t_n)\big) & \kern -.25cm \geq B_k(x_k)-B_{k+1}\big(x_{k+1}).
\ea\right.
\eeq
Hence, by virtue of condition \eqref{alkbek} the non-decreasing sequences $X_k(t_n)$ and $X_{k+1}(t_n)$ either are both bounded or both tend to $\infty$. This combined with $|X(t_n)|\to\infty$ thus implies that all the sequences $X_k(t_n)$ tend to $\infty$ as $n\to\infty$. As a consequence, since $u$ is separately non-decreasing, we get that for any $y\in\RR^d$,
\[
f(t_n)=u(X(t_n))\geq u(y)\quad\mbox{for any large enough $n$},
\]
which yields $\sup_{(\tau_-,\tau_+)}f=\sup_{\RR^d}u$. Similarly, we deduce from~\eqref{FkGk} that $\inf_{(\tau_-,\tau_+)}f=\inf_{\RR^d}u$. Therefore, since $f$ is increasing, we obtain the desired equality \eqref{rfu}.
\par
Now, fix a constant $c_u$ in the interval $(\inf_{\RR^d}u,\sup_{\RR^d}u)$. Then, for any $x\in\RR^d$, there exists a unique $\tau(x)\in (\tau_-,\tau_+)$ such that
\[
f(\tau(x))=u\big(X(\tau(x),x)\big)=c_u\in(\tex{\inf_{\RR^d}}u,\tex{\sup_{\RR^d}}u).
\]
Note that by virtue of the $C^2$-regularity of $u$, the flow $X(t,x)$ is a $C^1$-function (see, {\em e.g.}, \cite[Chap.~17.6]{HSD}) such that $\partial_t X$ is non-vanishing. Thus, the implicit functions theorem implies that $\tau$ belongs to $C^1(\RR^d)$.
\par
Then, the proof of the isotropic realizability of $\nabla u$ follows the same scheme that the proof of~\cite[Theorem~2.15]{BMT} with the time $\tau(x)$.
More precisely, by the semi-group property of the flow
\beq\label{sgX}
X(s,X(t,x))=X(s+t,x)\quad\mbox{for any $s,t$ close to 0},
\eeq
combined with the uniqueness of $\tau(x)$ we have for any $x\in\RR^d$,
\beq\label{tauX}
\tau\big(X(t,x)\big)=\tau(x)-t\quad\mbox{for any $t$ close to $0$}.
\eeq
Then, the $C^1$-function $\si$ defined by
\beq\label{sitau}
\sigma(x):=\exp\left(\int_0^{\tau(x)} \De u\big(X(s,x)\big)\,ds\right)\quad\mbox{for }x\in\RR^d,
\eeq
by \eqref{sgX} and \eqref{tauX} satisfies for any $t$ close to $0$, 
\beq\label{siX}
\sigma\big(X(t,x)\big)=\exp\left(\int_0^{\tau(x)-t} \De u\big(X(s+t,x)\big)\,ds\right)=\sigma(x)\exp\left(-\int_0^t \De u\big(X(s,x)\big)\,ds\right).
\eeq
Hence, differentiating the former equality with respect to $t$ and taking $t=0$, we get that for any $x\in\RR^d$,
\beq\label{siDu}
\nabla\sigma(x)\cdot\nabla u(x)=-\,\sigma(x)\,\De u(x)\quad\mbox{or equivalently}\quad{\rm div}\left(\sigma\nabla u\right)(x)=0.
\eeq
Therefore, $\nabla u$ is isotropically realizable in $\RR^d$ with the positive function $\si\in C^1(\RR^d)$.
\par\medskip\noindent
{\it Proof of $ii)$.} Let $x\in\RR^d$ be such that $u(x)>u(x^0)$. Let us prove that the range of the function $f$ defined by \eqref{f} contains the interval $(u(x^0),\sup_{\RR^d}u)$, {\em i.e.}
\beq\label{rfsupu}
\big(u(x^0),\tex{\sup_{\RR^d}u}\big)\subset \big\{f(t):t\in (\tau_-,\tau_+)\big\}.
\eeq
\par
First, note that
\beq\label{f'>0}
\forall\,t\in(\tau_-,\tau_+),\quad f'(t)=|\nabla(X(t,x))|^2>0.
\eeq
Indeed, if $f'(t_0)=0$ for some $t_0$, then $X(t_0,x)=x^0=X(t_0,x^0)$. Hence, by the uniqueness of the Cauchy-Lipschitz theorem, $X(t,x)=x^0$ for any $t\in(\tau_-,\tau_+)$ and $x=X(0,x)=x^0$, a contradiction.
\par
The inequality \eqref{f'>0} combined with the first argument of case $i)$ implies that the flow $X(t)$ is not bounded in the neighborhood of $\tau_+$. Thus, as in the case $i)$ with \eqref{FkGk} we get that $\sup_{(\tau_-,\tau_+)}f=\sup_{\RR^d}u$. Moreover, if the flow $X(t)$ is not bounded in the neighborhood of $\tau_-$, then we obtain similarly that $\inf_{(\tau_-,\tau_+)}f=\inf_{\RR^d}u$. In this case the range of $f$ thus agrees with $(\inf_{\RR^d}u,\sup_{\RR^d}u)$, which implies \eqref{rfsupu}.
\par
It thus remains to study the case where the flow $X(t)$ is bounded in the neighborhood of $\tau_-$, which implies that $\tau_-=-\infty$. Moreover, the function $f'$ is not bounded by below by a positive constant in the neighborhood of $-\infty$, otherwise $f(t)=u(X(t))$ tends to $-\infty$ as $t\to-\infty$. Hence, there exists a decreasing sequence $t_n\leq 0$ which tends to $-\infty$ such that $X(t_n)$ tends to some point $\bar{x}$ and $f'(t_n)=|\nabla u(X(t_n))|^2$ tends to $0$ as $n\to\infty$. At the limit we get that $\nabla u(\bar{x})=0$, which implies that $\bar{x}=x^0$ and $\inf_{(-\infty,\tau_+)}f=u(x^0)$. Therefore, by the increase of $f$ we obtain that the range of $f$ is $(u(x^0),\sup_{\RR^d}u)$, which establishes \eqref{rfsupu}.
\par
Fix a constant $c_u$ in $(u(x^0),\sup_{\RR^d}u)$. Then, for any $x\in\RR^d$ such that $u(x)>u(x^0)$, there exists by \eqref{rfsupu} a unique time $\tau(x)\in (\tau_-,\tau_+)$ such that
\[
f(\tau(x))=u\big(X(\tau(x),x)\big)=c_u\in\big(u(x^0),\tex{\sup_{\RR^d}}u\big).
\]
\par
Finally, we prove the isotropic realizability of $\nabla u$ in the open set $\{u>u(x^0)\}$ following the argument between \eqref{tauX} and \eqref{siDu} with the time $\tau(x)$. The proof of the isotropic realizability of $\nabla u$ in the open set $\{u<u(x^0)\}$ is quite similar.
\par\medskip\noindent
{\it Proof of $iii)$.} Let $x\in\RR^d$ be such that $c_j<u(x)<c_{j+1}$ for some $j=0,\dots,n$.
Repeating the arguments of $i)$ and $ii)$ we have the following alternative satisfied by the function $f$ defined by~\eqref{f}:
\begin{itemize}
\item $X(t)$ is not bounded in the neighborhood of $\tau_+$ (resp. $\tau_-$), then $\sup_{(\tau_-,\tau_+)}f=\sup_{\RR^d}u$ (resp. $\inf_{(\tau_-,\tau_+)}f=\inf_{\RR^d}u$),
\item $X(t)$ is bounded in the neighborhood of $\tau_+$ (resp. $\tau_-$), then $\tau_+=\infty$ (resp. $\tau_-=-\infty$), and
$\sup_{(\tau_-,\tau_+)}f\geq c_{j+1}$ (resp. $\inf_{(\tau_-,\tau_+)}f\leq c_j$).
\par
Contrary to case $ii)$, here we have only an inequality since $c_{j+1}$ (respectively $c_j$) is the smallest (respectively largest) critical value which can be attained asymptotically by the function $f$.
\end{itemize}
Hence, we deduce that
\beq\label{rfcj}
(c_j,c_{j+1})\subset \big\{f(t):t\in (\tau_-,\tau_+)\big\}.
\eeq
Finally, we conclude as before by considering a constant $c_u\in (c_j,c_{j+1})$, the time $\tau(x)$ such that $u(X(\tau(x),x))=c_u$, and the conductivity $\sigma$ defined by \eqref{sitau} in the open set $\{c_j<u<c_{j+1}\}$.
\cqfd
\subsection{Isotropic realizability of a vector field in $\RR^d$}
In this section we consider the isotropic realizability of a vector field $b\in C^1(\RR^d)^d$. Consider the flow associated with the vector field $b$ defined by \eqref{bflow}.
In the sequel, $0\in (\tau_-(x),\tau_+(x))$ denotes the maximal interval on which the solution $X(\cdot,x)$ to \eqref{bflow} is defined.
\par
In the spirit of the former proof we have the following extension of Theorem~\ref{thm.irg}.
\begin{Thm}\label{thm.irb}
Let $b\in C^1(\RR^d)^d$ and let $I$ be a non-empty open interval of $\RR$. Assume that there exists a function $u\in C^1(\RR^d)$ such that for any $x\in \{u\in I\}$, the function $f_x:=u(X(\cdot,x))$ satisfies $f'_x>0$ in $(\tau_-(x),\tau_+(x))$ and
\beq\label{Ifx}
\forall\,x\in \{u\in I\},\quad I\subset\big\{f_x(t):t\in(\tau_-(x),\tau_+(x))\big\},
\eeq
where $\{u\in I\}$ denotes the inverse image of $I$ by $u$.
Then, the vector field $b$ is isotropically realizable in the open set $\{u\in I\}$ with a positive function $\sigma\in C^1(\{u\in I\})$.
\end{Thm}
\par\noindent
{\bf Proof of Theorem~\ref{thm.irb}.}
Fix a constant $c_I$ in the interval $I$, and let $x\in \{u\in I\}$. By \eqref{Ifx} and the increase of $f_x$, there exists a unique $\tau(x)\in(\tau_-(x),\tau_+(x))$ such that
\beq\label{utaucI}
u(X(\tau(x),x))=c_I.
\eeq
By the semi-group property \eqref{sgX} of the flow $X$ combined with the uniqueness of $\tau$ the equality \eqref{tauX} still holds true.
Moreover, by the implicit functions theorem $\tau$ belongs to $C^1(\{u\in I\})$. Therefore, following \eqref{sitau}, \eqref{siX}, \eqref{siDu} with $b$ instead of $\nabla u$, the $C^1$-function $\sigma$ defined by
\beq\label{sitaub}
\sigma(x):=\exp\left(\int_0^{\tau(x)} ({\rm div}\,b)\big(X(s,x)\big)\,ds\right)\quad\mbox{for }x\in \{u\in I\},
\eeq
is solution to the equation ${\rm div}\,(\sigma b)=0$ in the open set $\{u\in I\}$.
\cqfd
\begin{Exa}\label{exa.irb}
Consider the gradient field $b=\nabla v$ in $\RR^2$ defined by
\[
v(x):={1\over 3}\,(x_1^3+x_2^3)\quad\mbox{for }x\in\RR^2.
\]
Then, the flow $X$ defined by \eqref{bflow} is given by
\[
X(t,x)=\left({x_1\over 1-t x_1},{x_2\over 1-t x_2}\right),\;\;t\in(\tau_-(x),\tau_+(x))=
\left\{\ba{cl}
(-\infty,{1\over\max(x_1,x_2)}) & \mbox{if }x_1>0,\,x_2>0
\\ \ecart
({1\over\min(x_1,x_2)},\infty) & \mbox{if }x_1<0,\,x_2<0
\\ \ecart
({1\over x_1},{1\over x_2}) & \mbox{if }x_1x_2<0
\\ \ecart
({1\over x_1},\pm\infty) & \mbox{if }\mp x_1>0,\,x_2=0
\\ \ecart
({1\over x_2},\pm\infty) & \mbox{if }x_1=0,\,\mp\,x_2>0
\\ \ecart
\RR & \mbox{if }x=(0,0).
\ea\right.
\]
It is clear that the function $v$ satisfies condition \eqref{dku>0} but not condition \eqref{disdku}.
\par
Define the function $u$ by $u(x):=x_1+x_2$ for $x\in\RR^2$. The function $f_x:=u(X(\cdot,x))$ satisfies for any $x\neq (0,0)$,
\[
\forall\,t\in (\tau_-(x),\tau_+(x)),\quad f'_x(t)=(\nabla u\cdot\nabla v)(X(t,x))=(X_1(t,x))^2+(X_2(t,x))^2>0,
\]
and
\[
\big\{f_x(t):t\in(\tau_-(x),\tau_+(x))\big\}=\left\{\ba{cl}
(0,\infty) & \mbox{if }x_1>0,\,x_2>0
\\ \ecart
(-\infty,0) & \mbox{if }x_1<0,\,x_2<0
\\ \ecart
\RR & \mbox{if }x_1x_2<0
\\ \ecart
(0,\pm\infty) & \mbox{if }\pm x_1>0,\,x_2=0
\\ \ecart
(0,\pm\infty) & \mbox{if }x_1=0,\,\pm\,x_2>0
\\ \ecart
\{0\} & \mbox{if }x=(0,0).
\ea\right.
\]
Hence, the function $u$ satisfies the conditions of Theorem~\ref{thm.irb} with $I=(0,\infty)$ and $I=(-\infty,0)$.
Define $c_I:=\pm 1$ if $I:=(0,\pm\infty)$, and let $x\in I$.
Moreover, it is easy to check that the solution $\tau(x)$ of \eqref{utaucI} is given by
\[
\tau(x)={x_1+x_2-2x_1x_2-\sqrt{(x_1-x_2)^2+4x_1^2x_2^2}\over 2x_1x_2}\quad\mbox{for }x_1+x_2\neq 0.
\]
Therefore, using formula \eqref{sitaub} we obtain that the gradient field $b=\nabla v$ is isotropically realizable in the open set $\{x_1+x_2\neq 0\}$ with the conductivity $\si\in C^1(\{x_1+x_2\neq 0\})$ defined by
\[
\sigma(x)={1\over\big(1-\tau(x)\,x_1\big)^2\big(1-\tau(x)\,x_2\big)^2}\quad\mbox{for }x_1+x_2\neq 0.
\]
\par
Note that $b=\nabla v$ is isotropically realizable in the open sets $\{x_1,x_2>0\}$ and $\{x_1,x_2<0\}$ with the simpler conductivity $x\mapsto (x_1x_2)^{-2}$. However, Theorem~\ref{thm.irb} here provides a suitable explicit conductivity in the two larger connected domains $\{x_1+x_2>0\}$ and $\{x_1+x_2<0\}$.
\end{Exa}
\section{Existence of a positive invariant measure in the torus under the gradient invertibility}
\subsection{Non-existence of a positive invariant measure}
In this section we will show that the isotropic realizability in $\RR^d$ of Theorem~\ref{thm.irg} under the positivity assumptions \eqref{dku>0} and \eqref{disdku} cannot be extended to the torus, namely the existence of a positive $Y_d$-periodic invariant measure.
In particular, note that any vector field $b\in C^0_\sharp(Y_d)^d$ satisfies conditions \eqref{dku>0} and \eqref{disdku} with $b$ in place of $\nabla u$, if there exists a constant $\al>0$ such that
\beq\label{b>0}
\forall\,k\in\{1,\dots,d\},\;\; b_k\geq \al\;\;\mbox{in }Y_d\quad\mbox{or}\quad \forall\,k\in\{1,\dots,d\},\;\;b_k\leq -\al\;\;\mbox{in }Y_d.
\eeq
\par
First, we have the following non-existence result if some component of $b$ changes sign.
\begin{Pro}\label{pro.neim}
Let $b$ a periodic vector field in $L^\infty_\sharp(Y_d)^d$.
Assume that for some $k=1,\dots,d$, say $k=1$ without loss of generality, there exist two measurable subsets $A_1$, $B_1$ of $[0,1]$ with positive Lebesgue measure, such that
\beq\label{AB}
\left\{\ba{lll}
\forall\,x_1\in A_1, & b_1(x_1,x')>0 & \mbox{a.e. }x'\in\RR^{d-1}
\\ \ecart
\forall\,x_1\in B_1, & b_1(x_1,x')<0 & \mbox{a.e. }x'\in\RR^{d-1}.
\ea\right.
\eeq
Then, the vector field $b$ has no positive invariant measure $\sigma\in L^\infty_\sharp(Y_d)$.
\end{Pro}
\noindent
{\bf Proof of Proposition~\ref{pro.neim}.} Assume by contradiction that there exists some positive function $\si\in L^\infty_\sharp(Y_d)^d$ such that ${\rm div}(\sigma b)=0$ in $\RR^d$, or equivalently in the torus sense
\[
\forall\,\ph\in H^1_\sharp(Y_d),\quad\int_{Y_d}\sigma(x)b(x)\cdot\nabla\ph(x)\,dx=0.
\]
If the function $\ph$ only depends on the variable $x_1$, the former equation leads us to
\[
\forall\,\ph\in H^1_\sharp(0,1),\quad\int_0^1\left(\int_{Y_{d-1}}\sigma(x_1,x')\,b_1(x_1,x')\,dx'\right)\ph'(x_1)\,dx=0,
\]
which implies the existence of a constant $c\in\RR$ such that
\[
\int_{Y_{d-1}}\sigma(x_1,x')\,b_1(x_1,x')\,dx'=c\quad\mbox{a.e. }x_1\in[0,1].
\]
Then, by assumption \eqref{AB} combined with the Fubini theorem we get that
\[
\left\{\ba{l}
\dis \int_{A_1\times Y_{d-1}}\underbrace{\sigma(x_1,x')\,b_1(x_1,x')}_{>0}\,dx=c\,|A_1|>0
\\ \ecart
\dis \int_{B_1\times Y_{d-1}}\underbrace{\sigma(x_1,x')\,b_1(x_1,x')}_{<0}\,dx=c\,|B_1|<0,
\ea\right.
\]
which yields a contradiction.
\cqfd
\par\bs
However, the positivity property \eqref{b>0} satisfied by a vector field $b$ is not sufficient to ensure the existence of a positive periodic invariant measure as shows the following example.
\begin{Exa}\label{exa.noim}
Consider the $Y_2$-periodic continuous gradient field $b=\nabla u$ defined in $\RR^2$ by
\[
u(x):=\al x_1+\al x_2+2\cos(2\pi x_1)\cos(2\pi x_2)\;\;\mbox{for }x\in\RR^2,\quad\mbox{where }\al>4\pi.
\]
We have $\partial_{x_k}u\geq\al-4\pi>0$ in $Y_2$ for $k=1,2$, such that the gradient field $b=\nabla u$ satisfies condition~\eqref{b>0}.
On the other hand, make the change of function
\[
v(y):=u(x)=\cos(4\pi y_1)+2\al y_2+\cos(4\pi y_2),\quad\mbox{where}\quad \left\{\ba{l} x_1=y_1+y_2 \\ x_2=y_2-y_1. \ea\right.
\]
Then, the gradient field $\nabla_y v$ is still $Y_2$-periodic.
Moreover, by virtue of Proposition~\ref{pro.neim} $\nabla_y v$ has not a $Y_2$-periodic positive invariant measure, since $\partial_{y_1} v$ only depends on the variable $y_1$ and changes sign. Also note that the orthogonality of the change of variables $x=Py$, where $PP^T=2\,I_2$, preserves the isotropy.
Indeed, for any $\sigma\in L^\infty_\sharp(Y_2)$, and for any $\ph\in C^\infty_c(\RR^2)$ and $\psi(y):=\ph(x)$, we have $\nabla_y\psi(y)=P^T\nabla_x\ph(x)$, and
\[
\ba{ll}
\dis \int_{\RR^2}\sigma(Py)\nabla_y v(y)\cdot\nabla_y\psi(y)\,dy
& \dis =\int_{\RR^2}\sigma(x)PP^T\nabla u(x)\cdot\nabla\ph(x)\,|\det P|^{-1}\,dx
\\ \ecart
& \dis =\int_{\RR^2}\sigma(x)\nabla u(x)\cdot\nabla\ph(x)\,dx,
\ea
\]
which implies that
\[
{\rm div}_y\big(\sigma(Py)\nabla_y v\big)=0\;\;\mbox{in }\RR^2\;\;\Leftrightarrow\;\;{\rm div}(\sigma\nabla u)=0\;\;\mbox{in }\RR^2,
\]
where the positive function $y\mapsto \sigma(Py)$ belongs to $L^\infty_\sharp(Y_2)$.
Hence, the gradient field $\nabla u$ cannot have a $Y_2$-periodic positive invariant measure since $\nabla_y v$ has not one.
Therefore, the $Y_2$-periodic gradient field $\nabla u$ satisfies condition \eqref{b>0}, but has not a $Y_2$-periodic positive invariant measure.
\end{Exa}
\subsection{Criterium for the existence of a positive invariant measure}
In this section we will give a criterium on a regular $Y_d$-periodic vector field $b$ so that it has a positive $Y_d$-periodic invariant measure.
Let $b$ be a periodic vector field in $C^1_\sharp(Y_d)^d$, and consider the associated flow $X$ defined by~\eqref{bflow}.
\par
We have the following result.
\begin{Thm}\label{thm.eim}
Let $b\in C^0_\sharp(Y_d)^d$.
Then, the following assertions are equivalent:
\begin{itemize}
\item[$i)$] There exist a positive function $\sigma\in C^0_\sharp(Y_d)^d$ and a vector field $V=(v_1,\dots,v_d)\in C^1(\RR^d)^d$ with $DV\in C^0_\sharp(Y_d)^{d\times d}$, such that
\beq\label{bDv1}
b\cdot\nabla v_1=1\;\;\mbox{in }Y_d,
\eeq
\beq\label{bDvk}
\sigma b=\left\{\ba{ll}
R_\perp\nabla v_2 & \mbox{if }d=2
\\ \ecart
\nabla v_2\times\cdots\times\nabla v_d & \mbox{if }d\geq 3,
\ea\right.
\quad\mbox{in }Y_d.
\eeq
\item[$ii)$] There exist a vector field $W\in C^1(\RR^d)^d$ and a non-zero vector $\xi\in\RR^d$ such that
\beq\label{DW}
DW\in C^0_\sharp(Y_d)^{d\times d}\mbox{ with }\langle DW\rangle=I_d,\quad \det\,(DW)\neq 0\;\;\mbox{in }Y_d,\quad (DW)^Tb=\xi\;\;\mbox{in }Y_d.
\eeq
\end{itemize}
\end{Thm}
\begin{Rem}
\hfill\par\smallskip\noindent
1. The gradient invertibility \eqref{bDv1} may seem rather sharp. But Proposition~\ref{pro.bv1} below gives some general cases for which it holds true.
\par\medskip\noindent
2. In dimension $d=2$ due to the representation of divergence free functions as orthogonal gradients, condition \eqref{bDvk} is equivalent to the fact that $\sigma$ is a positive $Y_2$-periodic invariant measure. In higher dimension condition \eqref{bDvk} only implies the existence of a positive $Y_d$-periodic invariant measure, since a divergence free vector field in $\RR^d$ with $d\geq 3$, is not necessarily of the form \eqref{bDvk}.
\end{Rem}
\begin{Pro}\label{pro.bv1}
\hfill\par\smallskip\noindent
$i)$ Let $b\in C^0_\sharp(Y_d)^d$. Assume that 
\beq\label{bk>0}
\exists\,k\in\{1,\dots,d\},\quad b_k(x)=b_k(x_k)>0\;\;\mbox{for }x\in\RR^d,
\eeq
then equality \eqref{bDv1} holds true.
\par\medskip\noindent
$ii)$ Let $b\in C^1_\sharp(Y_d)^d$. Assume that for some $k\in\{1,\dots,d\}$, there exists a function $u\in C^1(\RR^d)$ such that $\nabla u$ is $Y_d$-periodic, $b\cdot \nabla u>0$ in $Y^d$, and the mapping
\beq\label{Xuk}
\Phi:(t,x)\longmapsto {D_x X(t,x)\nabla u(X(t,x))\over(b\cdot\nabla u)(X(t,x))}\;\;\mbox{is bounded and uniformly continuous in }\RR\times\RR^d.
\eeq
Then, condition \eqref{bDv1} still holds true.
\par\medskip\noindent
$iii)$ Let $v_1,\dots,v_d$ be $d\geq 2$ functions in $C^1(\RR^d)$ such that
\beq\label{Dvkep}
\forall\,k\in\{1,\dots,d\},\quad\nabla v_k\mbox{ is $Y_d$-periodic}\;\;\mbox{and}\;\;\|\nabla v_k-e_k\|_{L^\infty(Y_d)^d}<\ep.
\eeq
Then, for any small enough $\ep>0$, the vector field $b\in C^0_\sharp(Y_d)^d$ defined by
\beq\label{bDvksidet}
\sigma:=\det\left(\nabla v_1,\nabla v_2,\dots,\nabla v_d\right)>0\quad\mbox{and}\quad
\sigma b:=\left\{\ba{ll}
R_\perp\nabla v_2 & \mbox{if }d=2
\\ \ecart
\nabla v_2\times\cdots\times\nabla v_d & \mbox{if }d>2,
\ea\right.
\eeq
satisfies condition \eqref{bDv1}.
\end{Pro}
\begin{Rem}
\hfill\par\smallskip\noindent
1. Condition \eqref{bk>0} is a particular case of \eqref{Xuk}. Indeed, assuming \eqref{bk>0} and choosing $u(x):=x_k$ we have
\beq\label{F}
X_k(t,x)=F^{-1}\big(t+F(x_k)\big)\quad\mbox{where}\quad F(t):=\int_0^{t}b_k^{-1}(s)\,ds.
\eeq
It follows that for any $(t,x)\in\RR\times\RR^d$,
\[
{D_x X(t,x)\nabla u(X(t,x))\over(b\cdot\nabla u)(X(t,x))}={\nabla_x X_k(t,x)\over b_k(X(t,x))}={1\over b_k(X(t,x))}\,{F'(x_k)\,e_k\over F'(X_k(t,x))}
={e_k\over b_k(x_k)},
\]
which clearly satisfies \eqref{Xuk}.
\par\medskip\noindent
2. In condition \eqref{Dvkep} we may replace the canonical basis by any basis $(\xi^1,\dots,\xi^d)$ of $\RR^d$ such that $\det\,(\xi^1,\dots,\xi^d)>0$.
\end{Rem}
\par\noindent
{\bf Proof of Theorem~\ref{thm.eim}.} We prove the case $d\geq 3$. The case $d=2$ is quite similar.
\par\smallskip\noindent
$i)\Rightarrow ii)$ By \eqref{bDvk} we have
\beq\label{bDv1DV}
\sigma=\sigma b\cdot\nabla v_1=\det\,(\nabla v_1,\dots,\nabla v_d)=\det\,(DV),
\eeq
which by the quasi-affinity of the cofactors (see, {\em e.g.}, \cite[Sec.~4.3.2]{Dac}) implies that
\[
\det\big(\langle DV\rangle\big)=\big\langle\!\det\,(DV)\big\rangle=\langle\sigma\rangle>0.
\]
Hence, the matrix $\langle DV\rangle$ is invertible, so that we may define the matrix $M$ of $\RR^{d\times d}$ by
\beq\label{M}
M:=\langle DV\rangle^{-1}.
\eeq
Let $W$ be the vector field defined by
\beq\label{W}
W:=M^TV \in C^1(\RR_d)^d,
\eeq
and let $\xi$ be the vector defined by
\beq\label{xi}
\xi:=M^Te_1.
\eeq
Then, by \eqref{bDv1}, \eqref{bDvk} and \eqref{xi} we get that
\beq\label{DWbxi}
(DW)^T b=M^T(DV)^T b=M^T\begin{pmatrix} b\cdot\nabla v_1 \\ b\cdot\nabla v_2 \\ \vdots \\ b\cdot\nabla v_d\end{pmatrix}=\xi.
\eeq
Moreover, by \eqref{M} and \eqref{W} we have $\langle DW\rangle=I_d$, and by \eqref{bDv1DV} we obtain that
\[
\det\,(DW)=\det\,(M)\,\det\,(DV)=\det\,(M)\,\sigma\neq 0\;\;\mbox{ in }Y_d.
\]
Therefore, the function $W$ satisfies the desired condition \eqref{DW}.
\par\medskip\noindent
$ii)\Rightarrow i)$ Let $W$ be a vector field satisfying \eqref{DW}. Consider an invertible matrix $M\in\RR^{d\times d}$ such that equation \eqref{xi} holds true, and define the vector field $V$ by~\eqref{W}. Then, we have the equalities \eqref{DWbxi} which combined with \eqref{xi} yield
\beq\label{bDv1vk}
b\cdot\nabla v_1=1\;\;\mbox{and}\;\; b\cdot\nabla v_2=\cdots=b\cdot\nabla v_d=0\quad\mbox{in }Y_d,
\eeq
which implies in particular \eqref{bDv1}.
Moreover, we have
\[
\det\,(DW)=\det\,(M)\,\det\,(DV)=\det\,(M)\,\nabla v_1\cdot(\nabla v_2\times\cdots\times\nabla v_d)\neq 0\quad\mbox{in }Y_d.
\]
Therefore, using a continuity argument and up to change $v_2$ in $-\,v_2$, the orthogonality conditions of \eqref{bDv1vk} imply the existence of a positive function $\sigma\in C^0_\sharp(Y_d)$ such that
condition \eqref{bDvk} holds true, which concludes the proof.
\cqfd
\par\bs\noindent
{\bf Proof of Proposition~\ref{pro.bv1}.}
\par\smallskip\noindent
{\it Proof of $i)$.} Assume that \eqref{bk>0} holds true for some $k\in\{1,\dots,d\}$, and let $v_1$ be the function defined by
\[
v_1(x):=\int_0^{x_k}b_k^{-1}(s)\,ds\quad\mbox{for }x\in\RR^d.
\]
Therefore, $\nabla v_1=b_k^{-1}e_k$ is $Y_d$-periodic, and $b\cdot\nabla v_1=1$ in $Y_d$.
\par\medskip\noindent
{\it Proof of $ii)$.} Let $x\in\RR^d$. Define the function $f$ by $f(t):=u(X(t,x))$ for $t\in\RR$. There exists a constant $c>0$ such that
\[
\forall\,t\in\RR,\quad f'(t)=(b\cdot\nabla u)(X(t,x))\geq c,
\]
which implies that the range of $f$ is $\RR$. Hence, there exists a unique $\tau(x)\in\RR$ such that $f(\tau(x))=0$, and by the implicit functions theorem $\tau$ belongs to $C^1(\RR)$. Moreover, by the semi-group property of the flow combined with the uniqueness of $\tau(x)$, we have
\[
\forall\,t\in\RR,\quad \tau(X(t,x))=\tau(x)-t.
\]
On the one hand, taking the derivative with respect to $t$ and choosing $t=0$, we get that
\beq\label{bDtau}
b\cdot\nabla\tau=-1\;\;\mbox{in }\RR^d.
\eeq
On the other hand, differentiating with respect to $x$ the equality $u(X(\tau(x),x))=0$, we get that
\[
(b\cdot\nabla u)(X(\tau(x),x))\nabla\tau(x)+D_x X(\tau(x),x)\nabla u(X(\tau(x),x))=0,
\]
or equivalently,
\[
\nabla\tau(x)=-\,{D_x X(\tau(x),x)\nabla u(X(\tau(x),x))\over (b\cdot\nabla u)(X(\tau(x),x))}=-\,\Phi(\tau(x),x).
\]
This combined with \eqref{Xuk} implies that $\nabla\tau$ is bounded and uniformly continuous in $\RR^d$.
Hence, by the Ascoli theorem the average of gradient functions
\[
x\longmapsto {-1\over(2n+1)^d}\kern -.2cm\sum_{\kappa\in\ZZ^d\cap[-n,n]^d}\kern -.4cm\nabla\tau(x+\kappa)
\]
converges uniformly, up to a subsequence of $n$, to some continuous gradient $\nabla v_1$ in any compact set of $\RR^d$.
The function $\nabla v_1$ is clearly $Y_d$-periodic, and equality \eqref{bDtau} implies~\eqref{bDv1}.
\par\medskip\noindent
{\it Proof of $iii)$.} Condition \eqref{Dvkep} implies that
\[
\sigma=\det\left(e_1+\nabla v_1-e_1,\dots,e_d+\nabla v_d-e_d\right)=1+O(\ep),
\]
so that $\sigma$ is positive when $\ep$ is small enough.
Then, the vector field $b$ defined by \eqref{bDvksidet} satisfies
\[
b\cdot\nabla v_1=\left\{\ba{ll}
\sigma^{-1} \nabla v_1\cdot R_\perp\nabla v_2 & \mbox{if }d=2
\\ \ecart
\sigma^{-1} \nabla v_1\cdot \left(\nabla v_2\times\cdots\times\nabla v_d\right) & \mbox{if }d>2
\ea\right\}
=\sigma^{-1}\det\left(\nabla v_1,\dots,\nabla v_d\right)=1\quad\mbox{in }Y_d,
\]
which concludes the proof.
\cqfd
\section{Applications}
\subsection{Asymptotics of the flow}
There exists an interesting by-product of Theorem~\ref{thm.eim} in terms of the asymptotics of the flow $X$ defined by \eqref{bflow}, which gives an alternative approach to the classical ergodic approach.
\par
We have the following result.
\begin{Cor}\label{cor.asyX}
Let $b\in C^1_\sharp(Y_d)^d$ be a vector field such that conditions \eqref{bDv1} and \eqref{bDvk} hold true with functions $v_k\in C^2(\RR^d)$.
Then, there exists a vector $\xi\in\RR^d$ such that the flow $X$ defined by \eqref{bflow} satisfies
\beq\label{asyX}
\forall\,x\in\RR^d,\quad\lim_{|t|\to\infty}{X(t,x)\over t}=\xi={\langle \sigma b\rangle\over\langle\sigma\rangle}
\eeq
Moreover, if there exists a non-zero vector $\la\in\RR^d$ such that $b\cdot\la=0$ in $Y_d$, and if either $\sigma b$ is not constant in dimension $d=2$ or $\sigma b$ is not of the form $\la\times\nabla w$ in dimension $d=3$, then the flow $X$ is not ergodic.
\end{Cor}
\begin{Rem}\label{rem.asyX}
\hfill\par\smallskip\noindent
1. Condition \eqref{bDvk} implies the existence of a positive $Y_d$-periodic invariant measure for the flow $X$ associated with the vector field $b\in C^1_\sharp(Y_d)^d$.
Hence, by virtue of the Birkhoff ergodic theorem (see, {\em e.g.}, \cite[Chap.~II, \S 5]{ReSi})
\beq\label{b*}
b^*(x):=\lim_{|t|\to\infty}{X(t,x)\over t}\quad\mbox{exists for a.e. }x\in\RR^d,
\eeq
with respect to the Lebesgue measure, and $b^*$ is invariant by the flow $X$, {\em i.e.}
\beq\label{b*inv}
b^*(X(t,x))=b^*(x)\quad\forall\,t\in\RR,\ \mbox{a.e. }x\in\RR^d.
\eeq
Under the extra assumption \eqref{bDv1} Corollary~\ref{cor.asyX} shows that limit \eqref{b*} holds actually for every $x\in\RR^d$. Moreover, the limit $b^*$ turns out to be the constant $\xi$ of \eqref{DW}, while the flow $X$ is not in general ergodic as shown in Example~\ref{exa.nonergo}.
Therefore, to prove \eqref{asyX} we need to use a different approach from the classical ergodic approach.
\par\medskip\noindent
2. Formula \eqref{F} shows that in dimension $d=1$, under condition~\eqref{bDv1} or equivalently assuming that $b$ is a non-vanishing function in $C^1_\sharp(\RR)$, the flow $X$ is ergodic and asymptotics \eqref{asyX} holds with $\sigma=b^{-1}$ and the harmonic mean $\xi=\langle b^{-1}\rangle^{-1}$.
\par\medskip\noindent
3. In the particular case of dimension two, assume that the vector field $b\in C^1_\sharp(Y_2)^2$ has a positive measure invariant $\sigma\in C^0_{\sharp}(Y_2)$ and does not vanish in $\RR^2$. Then, using the Kolmogorov theorem (see \cite[Lect.~11]{Sin} or \cite[Section~2]{Tas}) Tassa~\cite[Section~3]{Tas} obtained the following asymptotics
\[
\lim_{t\to\pm\infty}{X(t,x)\over t}=a^*(e_1+\gamma\,e_2)\quad\mbox{for any }x\in\RR^d,
\quad\mbox{where }\gamma:={\langle \sigma b_2\rangle\over\langle\sigma b_1\rangle},
\]
with the alternative according to the so-called rotation number $\ga$:
\begin{itemize}
\item If $\gamma\notin\QQ$, or equivalently the flow $X$ is ergodic, we have
\[
a^*={\langle \sigma b_1\rangle\over\langle\sigma\rangle}\quad\mbox{and}
\quad a^*(e_1+\gamma\,e_2)={\langle \sigma b\rangle\over\langle\sigma\rangle}.
\]
\item If $\gamma\in\QQ$, we have in general
\[
a^*(e_1+\gamma\,e_2)\neq {\langle \sigma b\rangle\over\langle\sigma\rangle}.
\]
In this case the gradient invertibility \eqref{bDv1} cannot hold.
\end{itemize}
\par
In view of the two points above, condition \eqref{bDv1} gives the same asymptotics~\eqref{asyX} than in the ergodicity setting, but does not imply the ergodicity of the flow. Moreover, the loss of condition~\eqref{bDv1} does not imply the loss of the ergodicity assumption.
Therefore, condition \eqref{bDv1} can be regarded as a substitute for the classical ergodicity assumption, since it induces a new and different regime for  getting \eqref{asyX}.
\par\medskip\noindent
4. Peirone \cite[Theorem~3.1]{Pei} proved the asymptotics \eqref{asyX} everywhere in~$\RR^2$ under the sole condition that the vector field $b\in C^1_\sharp(Y_2)^2$ is non-vanishing in $\RR^2$, using to this end the Birkhoff ergodic theorem combined with the Poincar\'e-Bendixson theorem (see, {\em e.g.}, \cite[Sec.~10.5]{HSD}). Moreover, he provided \cite[Lemma~4.6]{Pei} an example of a non-vanishing vector field $b$ in $\RR^3$ such that the asymptotics \eqref{asyX} does not hold at some point.
\par
Therefore, Corollary~\ref{cor.asyX} gives an alternative approach for proving~\eqref{asyX} in dimension two in the absence of ergodicity assumption, and in higher dimension condition \eqref{bDvk} gives a large class of vector fields $b$ such that \eqref{asyX} holds true everywhere in~$\RR^d$.
\end{Rem}
\begin{Exa}\label{exa.nonergo}
\hfill\par\smallskip\noindent
1. Let $v_1\in C^1(\RR)$ be a function such that $v_1'$ is positive and $1$-periodic, and let $v\in C^1(\RR)$ be a positive, $1$-periodic and non-constant function.
Define the vector field $b\in C^1_\sharp(Y_2)^2$ and the positive function $\sigma\in C^0_\sharp(Y_2)$ by
\[
b(x):={1\over v_1'(x_2)}\,e_2\;\;\mbox{and}\;\;\sigma(x):=v(x_1)\, v_1'(x_2)\quad\mbox{for }x\in\RR^2.
\]
We have $b\cdot e_1=0$ in $Y_2$, $\si b=v(x_1)\,e_2$ is non-constant and divergence free, and condition~\eqref{bDv1} holds.
Therefore, by virtue of Corollary~\ref{cor.asyX} the flow $X$ defined by \eqref{bflow} is not ergodic and satisfies for any $x\in\RR^2$,
\[
\lim_{|t|\to\infty}{X(t,x)\over t}={\langle \sigma b\rangle\over\langle\sigma\rangle}={\langle v(x_1)\,e_2\rangle\over\langle v(x_1)v'_1(x_2)\rangle}={\,e_2\over v_1(1)-v_1(0)}.
\]
\par\medskip\noindent
2. Let $v_1,v_2,v_3\in C^1(\RR^3)$ be $3$ functions such that condition \eqref{Dvkep} holds true with small enough $\ep>0$, and that $v_2(x)=v_2(x_2)$ has a $1$-periodic and non-constant derivative.
Then, define the positive function $\sigma\in C^0_\sharp(Y_3)$ and the vector field $b\in C^1_\sharp(Y_3)^3$ by
\[
\sigma:=\det\left(\nabla v_1,\nabla v_2,\nabla v_3\right)=\nabla v_1\cdot(\nabla v_2\times\nabla v_3)
\;\;\mbox{and}\;\; \sigma b=\nabla v_2\times\nabla v_3=v'_2(x_2)\,e_2\times\nabla v_3\quad\mbox{in }Y_3,
\]
so that $b\cdot\nabla v_1=1$ and $b\cdot e_2=0$ in $Y_3$.
\par
On the other hand, for any $v\in C^1(\RR^3)$ with $\nabla v$ $Y_3$-periodic, the functions $\sigma b$ and $e_2\times\nabla v$ cannot agree. Otherwise, we have
\[
v'_2(x_2)\big(\partial_{x_3}v_3(x)\,e_1-\partial_{x_1}v_3(x)\,e_3\big)=\partial_{x_3}v(x)\,e_1-\partial_{x_1}v(x)\,e_3\quad\mbox{for }x\in\RR^3,
\]
which implies that there exists a function $w\in C^0(\RR)$ such that
\[
v(x)=v'_2(x_2)\,v_3(x)+w(x_2)\quad\mbox{for }x\in\RR^3.
\]
Hence, by the $Y_3$-periodicity of $\nabla v_3$ and $\nabla v$ it follows that
\[
v(x_1,x_2,x_3+1)-v(x)=\langle\partial_{x_3}v\rangle=v'_2(x_2)\,\big(v_3(x_1,x_2,x_3+1)-v_3(x)\big)=v'_2(x_2)\,\langle\partial_{x_3}v_3\rangle
\quad\mbox{for }x\in\RR^3,
\]
a contradiction since $v'_2$ is not constant and $\langle\partial_{x_3}v_3\rangle$ is close to $1$.
\par
Therefore, by virtue of Corollary~\ref{cor.asyX} the flow $X$ defined by \eqref{bflow} is not ergodic, and by the quasi-affinity of the cofactors satisfies for $x\in\RR^3$,
\[
\lim_{|t|\to\infty}{X(t,x)\over t}={\langle \sigma b\rangle\over\langle\sigma\rangle}
={\langle\nabla v_2\rangle\times \langle\nabla v_3\rangle\over\langle\nabla v_1\rangle\cdot\left(\langle\nabla v_2\rangle\times \langle\nabla v_3\rangle\right)}
={\langle\partial_{x_3}v_3\rangle\,e_1-\langle\partial_{x_1}v_3\rangle\,e_3\over \langle\partial_{x_1}v_1\rangle\langle\partial_{x_3}v_3\rangle
-\langle\partial_{x_3}v_1\rangle\langle\partial_{x_1}v_3\rangle}.
\]
\end{Exa}
\noindent
{\bf Proof of Corollary~\ref{cor.asyX}.}
By virtue of Theorem~\ref{thm.eim} there exist a function $W\in C^1(\RR)^d$ and a non-zero vector $\xi\in\RR^d$ satisfying \eqref{DW}.
Define the function $W_\sharp\in C^1(\RR_d)^d$ by
\beq\label{Wd}
W_\sharp(x):=x-W(x)\quad\mbox{for }x\in\RR^d.
\eeq
Note that the vector field $W_\sharp$ is $Y_d$-periodic since $\langle DW\rangle=I_d$.
Then, using that
\[
b=(DW)^Tb+(DW_\sharp)^Tb=\xi+(DW_\sharp)^Tb,
\]
we have
\beq\label{asyXxiWd}
\ba{ll}
\dis X(t,x)=x+\int_0^{t} b(X(s,x))\,ds & \dis =x+t\,\xi+\int_0^{t} \big((DW_\sharp)^Tb\big)(X(s,x))\,ds
\\ \ecart
& \dis =x+t\,\xi+\int_0^{t} {\partial\over\partial s}\big(W_\sharp(X(s,x))\big)\,ds
\\ \ecart
& \dis =x+t\,\xi+W_\sharp(X(t,x))-W_\sharp(x).
\ea
\eeq
Since the function $W_\sharp\in C^1_\sharp(Y_d)^d$ is bounded in $\RR^d$, the former equality implies limit \eqref{asyX}.
\par\noindent
Moreover, by \eqref{DW} and the periodic div-curl lemma we have for any $\la\in\RR^d$ and $w_\la:=W\la$,
\[
\langle\sigma\rangle\,\xi\cdot\la=\langle\sigma\,\xi\cdot\la\rangle=\langle \sigma b\cdot\nabla w_\la\rangle=\langle\sigma b\rangle\cdot\langle \nabla w_\la\rangle
=\langle\sigma b\rangle\cdot\la,
\]
which yields the second equality of \eqref{asyX}.
\par\medskip
Assume that there exists a non-zero vector $\la\in\RR^d$ such that $b\cdot\la=0$ in $Y_d$, and that either $\sigma b$ is not constant in dimension $d=2$ or $\sigma b$ is not of the form $\la\times\nabla w$ in dimension $d=3$.
Then, using the quasi-affinity of the determinant and \eqref{cropro} we have
\[
\det\big(\la,\langle\nabla v_2\rangle,\dots,\langle\nabla v_d\rangle\big)
=\big\langle\la\cdot(\nabla v_2\times\cdots\times\nabla v_d)\big\rangle=\langle\sigma b\cdot\la\rangle=0,
\]
which implies the existence of a non-zero vector $(\al_1,\al_2,\dots,\al_d)\in\RR^d$ such that
\[
\al_1\,\la+\al_2\,\langle\nabla v_2\rangle+\cdots+\al_d\,\langle\nabla v_d\rangle=0.
\]
In view of \eqref{bDvk} it follows that the function $v$ defined by
\[
v(x):=\al_1\,\la\cdot x+\al_2\,v_2(x)+\cdots+\al_d\,v_d(x)\quad\mbox{for }x\in\RR^d,
\]
is in $C^1_\sharp(Y_d)$, and satisfies the equality $b\cdot\nabla v=0$ in~$Y_d$, or equivalently \eqref{b*inv} with $b^*=\nabla v$.
Moreover, due to $\la\neq 0$ one of the coefficients $\al_i$ for some $i\geq 2$, is not null, say $\al_2\neq 0$ without loss of generality.
\par
Now, assume that the function $v$ is constant. Then, when $d=2$ we have by~\eqref{bDvk}
\[
0=\nabla v=\al_1\la+\al_2\nabla v_2=\al_1\la-\al_2\,R_\perp\sigma b,
\]
a contradiction with the assumption that $\sigma b$ is not constant. When $d=3$ we have by \eqref{bDvk}
\[
0=\nabla v\times\nabla v_3 = \al_1\la\times\nabla v_3+\al_2\,\sigma b,
\]
again a contradiction with the assumption on $\sigma b$.
Hence, the function $v$ is a non-constant invariant periodic function for the flow $X$. Therefore, the flow $X$ is not ergodic (see \cite[Chap. II, \S 5]{ReSi}), which concludes the proof.
\cqfd
\subsection{Homogenization of a linear transport equation}
For $b\in C^1_\sharp(Y_d)^d$ and $u_0\in C^1(\RR^d)$, consider the following transport equation with an oscillating drift 
\beq\label{transpeq}
\left\{\ba{ll}
\dis {\partial u_\ep\over\partial t}-b(x/\ep)\cdot\nabla u_\ep=0 & \mbox{in }(0,\infty)\times\RR^d
\\ \ecart
u_\ep(0,x)=u_0(x) & \mbox{for }x\in\RR^d.
\ea\right.
\eeq
The homogenization of equation~\eqref{transpeq} was studied in the case of a two-dimensional divergence free vector field $b$ ({\em i.e.} with a constant invariant measure) through an ergodic approach by Brenier~\cite{Bre}, then by Hou and Xin~\cite{HoXi} with an oscillating initial datum which was specifically treated by a two-scale approach. Tassa~\cite{Tas} extended these results to any invariant measure in dimension two.
These results show that the ergodicity of the flow associated with $b$ leads us to a homogenized linear transport equation. In contrast, the loss of ergodicity implies that the limit of $u_\ep$ is not in general solution to a linear transport equation as Tartar \cite{Tar} showed. Here, using the non-ergodic approach of Corollary~\ref{cor.asyX} we obtain the following homogenization result in any dimension.
\begin{Cor}\label{cor.htranspeq}
Let $b\in C^1_\sharp(Y_d)^d$ be a vector field satisfying conditions \eqref{bDv1} and \eqref{bDvk} with functions $v_k\in C^2(\RR^d)$, and
let $u_0\in C^1(\RR^d)$.
Then, the solution $u_\ep$ in $L^1_{\rm loc}(\RR_+\times\RR^d)$ to the transport equation \eqref{transpeq} converges in $C^0_{\rm loc}(\RR_+\times\RR^d)$ to the function $u_0(x+t\xi)$ where the vector $\xi$ is given by \eqref{asyX}.
\end{Cor}
\par\noindent
{\bf Proof of Corollary~\ref{cor.htranspeq}.}
The flow $X_\ep$ associated with the oscillating vector field $b(x/\ep)$ is given by
\beq
X_\ep(t,x)=\ep\,X({t/\ep},{x/\ep})\quad\mbox{for }(t,x)\in\RR\times\RR^d,
\eeq
where $X$ is the flow associated with $b$.
Hence, due to the equality~\eqref{asyXxiWd} and the boundedness of the function $W_\sharp\in C^1_\sharp(Y_d)^d$ defined by \eqref{Wd} it follows that
\beq\label{estXep}
\dis X_\ep(t,x)=x+t\,\xi+\ep\,W_\sharp(X(t/\ep,x/\ep))-\ep\,W_\sharp(x/\ep)=x+t\,\xi+O(\ep).
\eeq
On the other hand, it is well known that the solution $u_\ep\in L^1_{\rm loc}(\RR_+\times\RR^d)$ to the transport equation \eqref{transpeq} is given by
\[
u_\ep(t,x)=u_0(X_\ep(t,x))\quad\mbox{for }(t,x)\in\RR_+\times\RR^d.
\]
Therefore, this combined with \eqref{estXep} and the continuity of $u_0$ implies that the sequence $u_\ep$ converges uniformly to $u_0(x+t\,\xi)$ in any compact set of $\RR_+\times\RR^d$.

\end{document}